\documentclass[12pt,a4paper]{amsart}
\usepackage{amsfonts}
\usepackage{amsthm}
\usepackage{amsmath}
\usepackage{amscd}
\usepackage[latin2]{inputenc}
\usepackage{t1enc}
\usepackage[mathscr]{eucal}
\usepackage{indentfirst}
\usepackage{graphicx}
\usepackage{graphics}
\usepackage{pict2e}
\usepackage{epic}
\numberwithin{equation}{section}
\usepackage[margin=3.0cm]{geometry}
\usepackage{epstopdf}

\theoremstyle{plain}
\newtheorem{Th}{Theorem}[section]
\newtheorem{Lemma}[Th]{Lemma}

 \theoremstyle{definition}

\newtheorem{Rem}[Th]{Remark}
\newtheorem{?}[Th]{Problem}

\begin{document}

\title{A Short Note on Asymptotic Enumeration of Contingency Tables with Non-Uniform Margins}

\author[Da Wu]{Da Wu}

\address{University of Pennsylvania \\ Department of Mathematics \\ David Rittenhouse Lab \\ 209 South 33rd Street \\ Philadelphia, PA, 10104-6395} 

\email{dawu@math.upenn.edu}



\begin{abstract} 
 In this short note, we compute the precise asymptotics for the number of contingency tables with non-uniform margins. More precisely, for parameter $n,\delta, B,C>0$, we consider the set of matrices whose first $[n^\delta]$ rows and columns have sum $[BCn]$ and the remaining $n$ rows and columns have sum $[Cn]$. We compute the precise asymptotics of the cardinality of this set when $B<B_c=1+\sqrt{1+1/C}$ using the maximal entropy principle introduced in \cite{Bar10}. The only contribution of this note is a detailed expansion of the determinant of quadratic forms in asymptotic formulas.  
\end{abstract}

\maketitle
\section{Introduction}
Let $\mathbf r=(r_1,\ldots, r_m)$ and $\mathbf c=(c_1,\ldots, c_n)$ be two positive integer vectors such that 
\begin{equation*}
	r_1+\ldots+r_m=c_1+\ldots+c_n=N.
\end{equation*}
Let $\mathscr M(\mathbf r, \mathbf c)$ be the set of $m\times n$ non-negative integer matrices with $i$-th row sum $r_i$ and $j$-th column sum $c_j$ for all $1\leq i\leq m, 1\leq j\leq n$. Suppose all the $r_i$ and $c_j$ depend on the dimension $m$ and $n$, one of the fundamental problems in Combinatorics is to provide the precise asymptotics of $\#\mathscr M(\mathbf r,\mathbf c)$ as $m,n\to \infty$. Recently, following the work by Pak and Lyu in \cite{DLP} we are interested in the case of non-uniform margin with two different values. More precisely, we consider the case when 
\begin{equation*}
	\widetilde{\mathbf r}=\widetilde{\mathbf c}=(\underbrace{[BCn],\ldots, [BCn]}_{\text{[$n^\delta]$ entries}},\underbrace{[Cn],\ldots, [Cn]}_{\text{$n$ entries}})\in \mathbb N^{[n^\delta]+n}
\end{equation*}
for parameters $B,C>0$ and $0\leq \delta<1$. Let 
\begin{equation*}
	\mathscr M_{n,\delta}(B,C):=\mathscr M(\widetilde{\mathbf r},\widetilde{\mathbf c}).
\end{equation*}  
We are interested in precise asymptotics of $\#\mathscr M_{n,\delta}(B,C)$ when $n\to \infty$. It is shown in \cite{DLP} that the typical table (defined in \cite{Bar10}) associated with $\widetilde{\mathbf r}$ and $\widetilde{\mathbf c}$ are uniform bounded in large $n$ limit when $B<B_c=1+\sqrt{1+1/C}$. In this case, we apply the maximal entropy method in \cite{Bar10} to compute the precise asymptotics of $\mathscr M_{n,\delta}(B,C)$. When $B>B_c$, since entries in top left corner will blow up in large $n$ limit, the precise asymptotics is not known. However, loose estimate of $\log \#\mathscr M_{n,\delta}(B,C)$ is known with error $O(n\log n+n^{2\delta})$, see the main theorem in \cite{LP}.  
\section{Precise Asymptotics of $\#\mathscr M_{n,\delta}(B,C)$ in sub-critical regime} 
In this section, we compute the precise asymptotic formula for $\#\mathscr M_{n,\delta}(B,C)$ when $0\leq\delta<1$ and $B<B_c=1+\sqrt{1+1/C}$ (subcritical case). The computation is based on Theorem $1.3$ in \cite{Bar10} and Lemma $5.1$ in \cite{DLP} , which will be restated below. 
\subsection{Review of Literature}
First, we recall the general asymptotic formula for $\#\mathscr M(\mathbf r,\mathbf c)$ when the all of the entries of typical table $Z$ are of same order. More detailed description can be found in \cite{Bar10} Section $1$. We say margins $\mathbf r=(r_1,\ldots,r_m)$ and $\mathbf c=(c_1,\ldots,c_n)$ are \textit{$\delta'$-smooth} if they satisfy the following two conditions:
\vskip0.1cm%
$(i)$: $m\geq\delta' n$ and $n\geq \delta' m$. Namely, dimensions of the matrix are of the same order asymptotically.
\vskip0.1cm%
$(ii)$: $\delta'\tau\leq z_{ij}\leq\tau$ for some $\tau$ such that $\tau\geq\delta'$
and all $1\leq i\leq m,1\leq j\leq n$. Namely, entries of typical table are of the same order asymptotically.
\vskip0.2cm%
Next, for typical table $Z=(z_{ij})$ associated with $\mathscr M(\mathbf r,\mathbf c)$, we define the quadratic form $q:\mathbb R^{m+n}\to\mathbb R$ as the following: 
\begin{equation*}
	q(s,t):=\frac{1}{2}\sum_{\substack{1\leq j\leq m \\ 1\leq k\leq n}}(z_{jk}^2+z_{jk})(s_j+t_k)^2,
\end{equation*}
where $s=(s_1,\ldots,s_m)$ and $t=(t_1,\ldots, t_n)$. Notice that the null space is spanned by the vector $\vec{u}=(\underbrace{1,\ldots,1}_{\text{$m$ entires}},\underbrace{-1,\ldots,-1}_{\text{$n$ entries}})$. Let $H=u^{\perp}\subseteq \mathbb R^{m+n}$ and $q|_{H}$ is a positive definite quadratic form and we can define its determinant $\det\left(q|_{H}\right)$ to be the product of non-zero eigenvalues of $q$. 
\par
We also define the polynomials $f,h:\mathbb R^{m+n}\to \mathbb R$ by
\begin{align}
	f(s,t):=\frac{1}{6}\sum_{\substack{1\leq j\leq m \\ 1\leq k\leq n}}z_{jk}(z_{jk}+1)(2z_{jk}+1)(s_j+t_k)^3
\end{align}
and 
\begin{align}
	h(s,t):=\frac{1}{24}\sum_{\substack{1\leq j\leq m \\ 1\leq k\leq n}}z_{jk}(z_{jk}+1)(6z_{jk}^2+6z_{jk}+1)(s_j+t_k)^4,
\end{align}
where $s=(s_1,\ldots,s_m)$ and $t=(t_1,\ldots,t_n)$. Consider the Gaussian probability measure on $H$ with density proportional to $e^{-q}$ and define 
\begin{equation*}
	\mu:=\mathbb E[f^2]\qquad\text{and}\qquad \nu:=\mathbb E[h].
\end{equation*}
Now, we can state the main theorem in \cite{Bar10} . 
\begin{Th}[\cite{Bar10}, Theorem $1.3$]
	Fix $0<\delta'<1$ and let $\mathbf r$ and $\mathbf c$ be $\delta'$-smooth margins and $Z=(z_{ij})$ be the associated typical table for $\mathscr M(\mathbf r,\mathbf c)$. Then 
	\begin{align}\label{general asymptotics of tables}
		\#\mathscr M(\mathbf r,\mathbf c)\asymp \frac{e^{g(Z)}\sqrt{m+n}}{(4\pi)^{(m+n-1)/2}\sqrt{\det\left(q|_{H}\right)}}\exp\left(-\frac{\mu}{2}+\nu\right)
	\end{align}
	as $m,n\to +\infty$. 
\end{Th}
\begin{Rem}
 There exists some positive constants $\gamma_1(\delta')$ and $\gamma_2(\delta')$ such that 
 \begin{equation*}
 	\gamma_1(\delta')\leq \exp\left(-\frac{\mu}{2}+\nu\right)\leq \gamma_2(\delta').
 \end{equation*}	
 Therefore, 
 \begin{equation}
 	\exp\left(-\frac{\mu}{2}+\nu\right)=O(1).
 \end{equation}
\end{Rem}
\begin{Rem}\label{Remark of coordinate basis change}
	Using the change of coordinate basis, 
\begin{align}
	\det\left(q|_H\right)=(m+n)\cdot 2^{1-m-n}\det Q,
\end{align}
where $Q=(q_{il})$ is the $(m+n-1)\times (m+n-1)$ symmetric matrix with 
\begin{align}\label{entries of $Q$}
\begin{split}
	q_{j,k+m}=q_{k+m, j} &=z_{jk}^2+z_{jk}\qquad\text{for}\  1\leq j\leq m, 1\leq k\leq n-1,\\
	q_{jj}&=r_j+\sum_{k=1}^n z_{jk}^2=\sum_{k=1}^n(z_{jk}+z_{jk}^2) \qquad\text{for}\ 1\leq j\leq m,\\
	q_{k+m,k+m} &=c_k+\sum_{j=1}^n z_{jk}^2=\sum_{j=1}^n (z_{jk}+z_{jk}^2) \qquad\text{for} \ 1\leq k\leq n-1.
\end{split}
\end{align}
Therefore, we can further simplify $(\ref{general asymptotics of tables})$ to 
\begin{align}\label{simplified version of asymptotics of tables}
	\#\mathscr M(\mathbf r,\mathbf c)\asymp \frac{e^{g(Z)}}{(2\pi)^{(m+n-1)/2}\sqrt{\det Q}}\exp\left(-\frac{\mu}{2}+\nu\right).
\end{align}
See \cite{Bar10} Section $1.4$ for a more detailed explanation. 
\end{Rem}
Next, we recall the key Lemma in \cite{DLP} regarding the asymptotics of entries of $Z=(z_{ij})$ associated with $\mathscr M_{n,\delta}(B,C)$. 
\begin{Lemma}[\cite{DLP}, Lemma $5.1$]\label{asymptotics of typical table}
	Fix $0\leq\delta<1$ and let $Z=(z_{ij})_{1\leq i,j\leq n+[n^\delta]}$ be the typical table of $\mathscr M_{n,\delta}(B,C)$. Let $B_c=1+\sqrt{1+1/C}$ and we have the following, \\
	$(i)$: If $B<B_c$, then 
	\begin{align}
		z_{11}=\frac{B^2(C+1)}{(B_c-B)(B_c+B-2)}+O(n^{\delta-1}),\qquad z_{1,n+1}=BC+O(n^{\delta-1}).
	\end{align}
	$(ii)$: If $B>B_c$, then 
	\begin{align}
		z_{n+1,n+1}=C+O(n^{\delta-1}),\qquad z_{1,n+1}=B_cC+O(n^{\delta-1}),\qquad n^{\delta-1}z_{11}=C(B-B_c)+O(n^{\delta-1}).
	\end{align}
\end{Lemma}
\begin{Rem}
	The behaviour of $z_{n+1,n+1}$ is more predictable. It is shown in \cite{Bar10} that 
	\begin{align}
		|z_{n+1,n+1}-C|=n^{\delta-1}z_{1,n+1}\leq BCn^{\delta-1}.
	\end{align}
	Hence, it is trivial that 
	\begin{align}
		z_{n+1,n+1}=C+O(n^{\delta-1}).
	\end{align}
\end{Rem}
\subsection{Computation of $\#\mathscr M_{n,\delta}(B,C)$}
Now, we go back to our setting of $\#\mathscr M_{n,\delta}(B,C)$. Recall $0\leq \delta<1$ and $B<B_c=1+\sqrt{1+1/C}$. First, notice that when $B<B_c$, all of entries of $Z=(z_{ij})$ have well-defined finite limits, and by symmetry
\begin{align}
\begin{split}
	e^{g(Z)} &=\prod_{1\leq i,j\leq n+[n^\delta]}\frac{(z_{ij}+1)^{z_{ij}+1}}{z_{ij}^{z_{ij}}}\\
	&=\left(\frac{(z_{11}+1)^{z_{11}+1}}{z_{11}^{z_{11}}}\right)^{[n^\delta]^2}\left(\frac{(z_{n+1,n+1}+1)^{z_{n+1,n+1}+1}}{z_{n+1,n+1}^{z_{n+1,n+1}}} \right)^{n^2}\left(\frac{(z_{1,n+1}+1)^{z_{1,n+1}+1}}{z_{1,n+1}^{z_{1,n+1}}}\right)^{2n[n^\delta]}
\end{split}
\end{align}
 Next, we compute the determinant of $Q$ in $(\ref{simplified version of asymptotics of tables})$. By $(\ref{entries of $Q$})$,  $Q$ has entries 
\begin{align}\label{Diagonal Entries 1}
	\begin{split}
		q_{jj} &=[BCn]+[n^\delta]\left(\frac{B^2(C+1)}{(B_c-B)(B_c+B-2)}+O(n^{\delta-1})\right)^2+n\left(BC+O(n^{\delta-1})\right)^2
	\end{split}
\end{align}
when $1\leq j\leq [n^\delta]$ and $[n^\delta]+n+1\leq j\leq 2[n^\delta]+n$ 
\begin{align}\label{Diagonal Entries 2}
	\begin{split}
		q_{jj} &=[Cn]+[n^\delta]\left(BC+O(n^{\delta-1})\right)^2+n(C+O(n^{\delta-1}))^2
	\end{split}
\end{align}
when $[n^\delta]+1\leq j\leq [n^\delta]+n$ and $2[n^\delta]+n+1\leq j\leq 2([n^\delta]+n)-1$. 
\begin{align}\label{off-diagonal 1}
	\begin{split}
		q_{ij}=q_{ji} &=\left(\frac{B^2(1+C)}{(B_c-B)(B_c+B-2)}+O(n^{\delta-1})\right)^2+\frac{B^2(1+C)}{(B_c-B)(B_c+B-2)}+O(n^{\delta-1})
	\end{split}
\end{align} 
when $1\leq i\leq [n^\delta]$ and $[n^\delta]+n+1\leq j\leq 2[n^\delta]+n$. 
\begin{align}\label{off-diagonal 2}
	\begin{split}
		q_{ij}=q_{ji}&=\left(BC+O(n^{\delta-1})\right)^2+BC+O(n^{\delta-1})
	\end{split}
\end{align}
when $1\leq i\leq [n^{\delta}]$, $2[n^\delta]+n+1\leq j\leq 2([n^{\delta}]+n)-1$ and when $[n^\delta]+1\leq i\leq [n^\delta]+n$, $[n^\delta]+n+1\leq j\leq 2[n^\delta]+n$.
\begin{align}\label{off-diagonal 3}
	\begin{split}
		q_{ij}=q_{ij} &=(C+O(n^{\delta-1}))^2+C+O(n^{\delta-1})
	\end{split}
\end{align}
when $[n^\delta]+1\leq i\leq [n^\delta]+n$ and $2[n^\delta]+n+1\leq j\leq 2([n^\delta]+n)-1$. 
\vskip0.1cm%
The rest of the entries are zero. Notice that all the off-diagonal entries have size $O(1)$ while all the entries on the diagonal has asymptotical order $n$. To compute the asymptotics of $\det Q$, we write $Q=A+E$ where $A=\text{diag}(q_{11},q_{22}, \ldots, q_{2([n^\delta]+n)-1,2([n^\delta]+n)-1})$ is the diagonal  matrix. By diagonal expansion of the determinant, 
\begin{align}
\begin{split}
	\det(Q) &=\det(A+E)\\
	&=\det(A)+S_1+S_2+\ldots+S_{2([n^\delta]+n-1)}+\det(E)
\end{split}
\end{align}
where 
\begin{align*}
	S_k=\sum_{1\leq i_1<\ldots<i_k\leq 2([n^\delta]+n)-1}\left(\prod_{r=1}^kq_{i_r,i_r}\right)\det\left(E_{i_1,\ldots,i_k}\right)
\end{align*}
$E_{i_1,\ldots,i_k}$ is the principle minor of order $2([n^\delta]+n)-1-k$ of $E$. Trivially, $\det(E)=0$, and   
\begin{align}
	\det(A)=\prod_{i=1}^{2([n^\delta]+n)-1}q_{ii}
\end{align}
Furthermore, $S_{2([n^\delta]+n-1)}=0$, and 
\begin{align*}
	&S_{2([n^\delta]+n)-3}\\
	&=\left(q_{11}\ldots q_{[n^\delta]+n-1,[n^\delta]+n-1} \right)\left(q_{[n^\delta]+n+2,[n^\delta]+n+2}\ldots q_{2[n^\delta]+2n-1,2[n^\delta]+2n-1} \right)\left[-(z_{1,n+1}^2+z_{1,n+1})^2 \right]\\
	&S_{2([n^\delta]+n)-4}\\
	&=0\\
	&S_{2([n^\delta]+n)-5}\\
	&=\left([n^\delta]\right)^2n(n-1)\left[(z_{11}^2+z_{11})(z_{n+1,n+1}^2+z_{n+1,n+1})-(z_{1,n+1}^2+z_{1,n+1})^2 \right]\\
	&S_{2([n^\delta]+n)-i}\\
	&=0\qquad\forall i\geq 6
\end{align*}
The above computation is based on the symmetry of typical table, i.e. 
\begin{align*}
	z_{ij}=z_{i'j'}\qquad\text{if $r_i=r_{i'}$ and $c_j=c_{j'}$}
\end{align*} 
Therefore, 
\begin{align}
	\begin{split}
		&\det(Q)\\
		&=\left(\prod_{i=1}^{2([n^\delta]+n)-1}q_{i}\right)-\left(q_{1}\ldots q_{[n^\delta]+n-1} \right)\left(q_{[n^\delta]+n+2}\ldots q_{2[n^\delta]+2n-1} \right)\left[(z_{1,n+1}^2+z_{1,n+1})^2 \right]\\
		&+\left([n^\delta]\right)^2n(n-1)\left[(z_{11}^2+z_{11})(z_{n+1,n+1}^2+z_{n+1,n+1})-(z_{1,n+1}^2+z_{1,n+1})^2 \right]
	\end{split}
\end{align}
where we write $q_i$ in place of $q_{ii}$. Finally, by (\ref{Diagonal Entries 1}), (\ref{Diagonal Entries 2}), (\ref{off-diagonal 1}), (\ref{off-diagonal 2}), (\ref{off-diagonal 3}), (\ref{simplified version of asymptotics of tables}) and Lemma \ref{asymptotics of typical table}, we get the precise asymptotics of $\#\mathscr M_{n,\delta}(B,C)$. 
\section{Left half of the $\mathscr M_{n,\delta}(B,C)$}
In this section, we compute the case when
\begin{align*}
	\mathbf r_1=(\underbrace{B_cCn,\ldots,B_cCn}_{\text{$[n^\delta]$ entries}},\underbrace{Cn-B_cCn^\delta,\ldots,Cn-B_cCn^\delta}_{\text{$n$ entires}} )\in \mathbb Z_{>0}^{n+n^\delta}
\end{align*}
and
\begin{align*}
	\mathbf c_1=(Cn,\ldots,Cn)\in \mathbb Z_{>0}^n
\end{align*}
By symmetry and margin conditions, $Z=(z_{ij})$ satisfies
\begin{align*}
\begin{cases}
	nz_{11}=B_c Cn\\
	n^\delta z_{11}+nz_{1,n+1}=Cn
\end{cases}
\end{align*}
which implies that 
\begin{align*}
	\begin{cases}
		z_{11}=B_cC\\
		z_{1,n+1}=C-z_{11}n^{\delta-1}=C-B_cC n^{\delta-1}
	\end{cases}
\end{align*}
Next, we compute the exact asymptotic formula of $\#\mathscr M(\mathbf r_1,\mathbf c_1)$. Recall the formula, 
\begin{align*}
	\frac{e^{g(Z)}}{(2\pi)^{(m+n-1)/2}\sqrt{\det Q}}\exp\left(-\frac{\mu}{2}+\nu\right)
\end{align*}
where $Q=(q_{ij})\in \mathbb R_{\geq 0}^{(2n+n^\delta-1)\times (2n+n^\delta-1)}$ has entries 
\begin{align*}
	q_{ii}=
	\begin{cases}
		B_c Cn+n\left(B_c C\right)^2\qquad & 1\leq i\leq [n^\delta]\\
		Cn-B_c Cn^\delta+n\left(C-B_c Cn^{\delta-1}\right)^2 \qquad & [n^\delta]+1\leq i\leq [n^\delta]+n\\
		Cn+n^\delta (B_c C)^2+(n-n^\delta)\left(C-B_c Cn^{\delta-1}\right)^2\qquad & [n^\delta]+n+1\leq i\leq 2n+[n^\delta]-1
	\end{cases}
\end{align*}
and 
\begin{align*}
	&q_{ij}=q_{ji}=\\
	&\begin{cases}
		B_c C+\left(B_c C\right)^2\qquad & 1\leq i\leq [n^\delta], n+[n^\delta]+1\leq j\leq 2n+[n^\delta]-1\\
		C-B_c Cn^{\delta-1}+\left(C-B_c Cn^{\delta-1}\right)^2\qquad & [n^\delta]+1\leq i\leq [n^\delta]+n, n+[n^\delta]+1\leq j\leq 2n+[n^\delta]-1
	\end{cases}
\end{align*}
The rest of the entries are $0$. We write $Q=A+E$ where
\begin{align*}
	A=\text{diag}\left(q_{11},\ldots,q_{2n+[n^\delta]-1,2n+[n^\delta]-1} \right)
\end{align*}
By diagonal expansion of the determinants, 
\begin{align}
\begin{split}
	\det(Q) &=\det(A+E)\\
	&=\det(A)+S_1+S_2+\ldots+S_{2n+[n^\delta]-2}+\det(E)
\end{split}
\end{align}
where 
\begin{align*}
	S_k=\sum_{1\leq i_1<\ldots<i_k\leq 2n+[n^\delta]-2}\left(\prod_{r=1}^kq_{i_r,i_r}\right)\det\left(E_{i_1,\ldots,i_k}\right)
\end{align*}
$E_{i_1,\ldots,i_k}$ is the principle minor of order $2n+[n^\delta]-2-k$ of $E$. It is not hard to see that 
\begin{align*}
	&S_{2n+[n^\delta]-2}=0\\
	&S_{2n+[n^\delta]-3}\\
	 &=-[n^\delta](n-1)\left\lbrace B_c Cn+n\left(B_c C\right)^2 \right\rbrace^{[n^\delta]-1}\left\lbrace Cn-B_c Cn^\delta+n\left(C-B_c Cn^{\delta-1}\right)^2\right\rbrace^n\times \\
	& \left\lbrace Cn+n^\delta (B_c C)^2+(n-n^\delta)\left(C-B_c Cn^{\delta-1}\right)^2 \right\rbrace^{n-2}\left(B_c C+(B_cC)^2\right)^2\\
	&-n(n-1)\left\lbrace B_c Cn+n\left(B_c C\right)^2 \right\rbrace^{[n^\delta]}\left\lbrace Cn-B_c Cn^\delta+n\left(C-B_c Cn^{\delta-1}\right)^2\right\rbrace^{n-1}\times \\
	& \left\lbrace Cn+n^\delta (B_c C)^2+(n-n^\delta)\left(C-B_c Cn^{\delta-1}\right)^2 \right\rbrace^{n-2} \left(C-B_c Cn^{\delta-1}+\left(C-B_c Cn^{\delta-1}\right)^2 \right)^2\\
	&S_{2n+[n^\delta]-i}=0
\end{align*} 
for all $i\geq 4$. Therefore, 
\begin{align*}
	&\det Q\\
	&=\det A+S_{2n+[n^\delta]-3}\\
	&=(q_{11})^{[n^\delta]}\left(q_{n+1,n+1}\right)^{n}\left(q_{2n+1,2n+1}\right)^{n-1}\\
	&-[n^\delta](n-1)(q_{11})^{[n^\delta]-1}\left(q_{n+1,n+1}\right)^{n}\left(q_{2n+1,2n+1}\right)^{n-2}\left(B_c C+(B_cC)^2\right)^2 \\
	&-n(n-1)(q_{11})^{[n^\delta]}\left(q_{n+1,n+1}\right)^{n-1}\left(q_{2n+1,2n+1}\right)^{n-2}\left(C-B_c Cn^{\delta-1}+\left(C-B_c Cn^{\delta-1}\right)^2 \right)^2    \\
	&=\left\lbrace B_c Cn+n\left(B_c C\right)^2\right\rbrace^{[n^\delta]}\left\lbrace Cn-B_c Cn^\delta+n\left(C-B_c Cn^{\delta-1}\right)^2\right\rbrace^n\times\\
	& \left\lbrace  Cn+n^\delta (B_c C)^2+(n-n^\delta)\left(C-B_c Cn^{\delta-1}\right)^2\right\rbrace^{n-1}\\
	&-[n^\delta](n-1)\left\lbrace B_c Cn+n\left(B_c C\right)^2 \right\rbrace^{[n^\delta]-1}\left\lbrace Cn-B_c Cn^\delta+n\left(C-B_c Cn^{\delta-1}\right)^2\right\rbrace^n\times \\
	& \left\lbrace Cn+n^\delta (B_c C)^2+(n-n^\delta)\left(C-B_c Cn^{\delta-1}\right)^2 \right\rbrace^{n-2}\left(B_c C+(B_cC)^2\right)^2\\
	&-n(n-1)\left\lbrace B_c Cn+n\left(B_c C\right)^2 \right\rbrace^{[n^\delta]}\left\lbrace Cn-B_c Cn^\delta+n\left(C-B_c Cn^{\delta-1}\right)^2\right\rbrace^{n-1}\times \\
	& \left\lbrace Cn+n^\delta (B_c C)^2+(n-n^\delta)\left(C-B_c Cn^{\delta-1}\right)^2 \right\rbrace^{n-2} \left(C-B_c Cn^{\delta-1}+\left(C-B_c Cn^{\delta-1}\right)^2 \right)^2.
\end{align*}
Plugging in (\ref{simplified version of asymptotics of tables}) and we are done.

\end{document}